\documentclass[12pt,a4paper]{article}
\usepackage{amsmath,amssymb,theorem}

\newcommand{\frA}{\mathfrak{A}}
\newcommand{\frB}{\mathfrak{B}}

\newcommand{\fra}{\mathfrak{a}}
\newcommand{\frb}{\mathfrak{b}}

\newcommand{\Z}{\mathbb{Z}}
\newcommand{\Q}{\mathbb{Q}}
\newcommand{\R}{\mathbb{R}}
\newcommand{\C}{\mathbb{C}}
\renewcommand{\H}{\mathbb{H}}
\newcommand{\eps}{\varepsilon}
\newcommand{\sumpr}{\sideset{}{'}\sum}
\newcommand{\N}{\mathbf{N}}
\newcommand{\dif}{\mathfrak{d}}
\newcommand{\D}{\mathbb{D}}
\newcommand{\absN}[1]{\bigl|\N(#1)\bigr|}
\newcommand{\dtimes}{d^\times\!}

\DeclareMathOperator{\re}{\mathrm{Re}}
\DeclareMathOperator{\im}{\mathrm{Im}}
\DeclareMathOperator{\Tr}{\mathrm{Tr}}
\DeclareMathOperator{\Res}{\mathrm{Res}}
\DeclareMathOperator{\CT}{\mathrm{CT}}

\numberwithin{equation}{subsection}

\theorembodyfont{\slshape}
\newtheorem{thm}{Theorem}[subsection]
\newtheorem{prop}[thm]{Proposition}
\newtheorem{cor}[thm]{Corollary}
\newtheorem{lem}[thm]{Lemma}
\newtheorem{defn}[thm]{Definition}
\newtheorem{rem}[thm]{Remark}

\newenvironment{proof}
	{\textbf{Proof.}}
	{\hfill\rule{0.5em}{2ex}\par\bigbreak}

\begin{document}
\title{Hecke's integral formula for quadratic extensions of a number field}
\author{Shuji Yamamoto}
\maketitle

\begin{abstract}
Let $K/F$ be a quadratic extension of number fields. 
After developing a theory of the Eisenstein series over $F$, 
we prove a formula which expresses a partial zeta function of $K$ 
as a certain integral of the Eisenstein series. 
As an application, we obtain a limit formula of Kronecker's type 
which relates the $0$-th Laurent coefficients at $s=1$ 
of zeta functions of $K$ and $F$. 
\end{abstract}

\section{Introduction}
Let $E(z,s)$ be the real analytic Eisenstein series defined by 
\[E(z,s)=\frac{1}{2}\sumpr_{m,n\in\Z}\frac{y^s}{\lvert mz+n\rvert^{2s}}
\qquad\bigl(y=\im z>0,\,\re(s)>1\bigr), \]
where $\sumpr$ means that the sum is taken except for $(m,n)=(0,0)$. 
Then it is classical that, when $z$ is an element of 
an imaginary quadratic field, $E(z,s)$ represents a zeta function 
of that field. More precisely, let $A$ be an ideal class of 
an imaginary quadratic field $K$, and $\frA$ an element of $A^{-1}$ 
of the form $\frA=\Z z+\Z$. We fix an embedding of $K$ into $\C$ and 
assume that $\im(z)>0$. Then the partial zeta function 
\[\zeta_K(s,A):=\sum_{\substack{\frB\in A\\\frB\subset O_K}}\N(\frB)^{-s}\]
can be written as 
\begin{equation}\label{eq:ImagQuad}
\zeta_K(s,A)=\frac{2}{w_K}\biggl(\frac{\sqrt{d_K}}{2}\biggr)^{-s}E(z,s),
\end{equation}
where $w_K$ and $d_K$ denote the number of roots of unity in $K$ and 
the absolute value of the discriminant of $K$, respectively. 

Hecke \cite{Hecke17} discovered an analogous formula 
for real quadratic fields. Now let $K$ be real quadratic, 
and $\frA=\Z z+\Z\in A^{-1}$ as before. Again we fix an embedding 
$K\hookrightarrow \R$ and denote the conjugate of $x\in K$ over $\Q$ 
by $x'$. Assume that $z'>z$. Then Hecke's integral formula is 
\begin{equation}\label{eq:RealQuad}
\zeta_K(s,A)=2d_K^{-s/2}\frac{\Gamma(s)}{\Gamma(s/2)^2}
\int_1^{\eps^2}E(z_t,s)\frac{dt}{t},
\end{equation}
where $\eps>1$ is the fundamental unit of $K$ and 
\[z_t:=\frac{t^{1/2}z+t^{-1/2}z'i}{t^{1/2}+t^{-1/2}i}. \]
For further discussions on this formula and related results, 
see Meyer \cite{Meyer57}, Siegel \cite{Siegel80} and 
a recent work of Manin \cite[\S2]{Manin04}. 

There are some results related to the formulas 
(\ref{eq:ImagQuad}) and (\ref{eq:RealQuad}). 
For example, when $F$ is a totally real field and $K$ is its CM extension, 
there is a formula relating the zeta functions of $K$ 
and the Eisenstein series over $F$ evaluated at certain CM-points, 
which is a generalization of (\ref{eq:ImagQuad}) 
(see Yoshida \cite{Yoshida03}). 
Another example was given by Konno \cite{Konno89}, who 
found a formula analogous to (\ref{eq:RealQuad}), 
which expresses a zeta function of $K$ as an integral 
of the Eisenstein series over $F$, when $F$ is imaginary quadratic 
and $K\supset F$ is absolutely biquadratic. 

In this paper, we consider the most general situation, i.e., 
an arbitrary quadratic extension $K/F$ of number fields. 
Let $T_{K/F}$ be the subgroup of $(K\otimes_\Q\R)^\times$ 
consisting of the elements $u$ such that $\N_{K/F}(u)=1$, 
and $U_{K/F}$ the intersection of $T_{K/F}$ and the unit group of $K$. 
Then our main result is the following: 
\begin{thm}[=Theorem \ref{thm:HIF}]\label{thm:IntroHIF}
If $A$ is a wide ideal class of $K$ and $\frA$ is an element of $A^{-1}$, 
then we have 
\[\xi_K(s,A)=\frac{1}{w_{K/F}}
\int_{T_{K/F}/U_{K/F}^2}\widehat{E}\bigl(\rho(\tilde{u}\frA),s\bigr)
\dtimes u. \]
Here $\xi_K(s,A)$ denotes the completed zeta function associated with $A$, 
and $\widehat{E}$ is the completed Eisenstein series over $F$ 
(precise definitions are given in \S2). 
\end{thm}
For the other notations used in the above theorem, see \ref{subsec:HIF}. 

\bigbreak

The contents of this paper is as follows. 
We develop a theory of the Eisenstein series over an arbitrary $F$ in \S2. 
After the definitions (\ref{subsec:Lattices} and \ref{subsec:Eisenstein}), 
we prove the functional equation (Theorem \ref{thm:FE}), 
the Fourier expansion (Theorem \ref{thm:FourierExp}) and 
the Kronecker limit formula (Theorem \ref{thm:EKLF}). 
Note that we consider the Eisenstein series as a function 
of a lattice in a certain vector space over $\R$, 
although there is a more traditional notion studied 
by Asai \cite{Asai70} and Jorgenson-Lang \cite{Jorgenson-Lang99}. 
(In fact, these two formulations are essentially equivalent. 
See Remark \ref{rem:Equivalent}.) 

In \S3, we prove our generalization of 
Hecke's integral formula (Theorem \ref{thm:HIF}). 
We also apply it to the Kronecker limit formula about 
the constant terms in the Laurent expansions 
of zeta functions at $s=1$ (Theorem \ref{thm:zetaKLF}). 
The result has a relative nature, 
in the sense that it compares zeta functions of $K$ and $F$.

\subsection{Notation}
As usual, the symbols $\Z$, $\Q$, $\R$ and $\C$ mean the ring of integers, 
rationals, real numbers and complex numbers, respectively. 
We also denote by $\H$ the quaternion division algebra of Hamilton: 
\[\H=\C\oplus\C j=\R\oplus\R i\oplus\R j\oplus\R k. \]
By $\C^\times_1$, we mean the group of complex numbers of absolute value $1$. 

\bigbreak

For an arbitrary number field $F$ (of finite degree), 
we use the following notations: 

$O_F$, $U_F$, $\dif_F$ and $d_F$ mean the ring of integers, 
the group of units in $F$, the different and 
the absolute value of the discriminant, respectively. 

We denote the set of infinite places of $F$ by $S_F$, 
and the subset of real (resp.\ complex) ones by $S_F^1$ (resp.\ $S_F^2$). 

For each $v\in S_F$, we denote the corresponding completion by $F_v$, 
and the embedding of $F$ into $F_v$ by $x\mapsto x_v$. 
The same notation is used to indicate the $v$-th component of 
an element $x\in F_\R:=F\otimes_\Q\R=\prod_{v\in S_F}F_v$. 
Moreover, we equip $F_v$ with the Lebesgue measure (resp.\ twice the 
Lebesgue measure) when $v$ is real (resp.\ complex). 

We denote the absolute norm (of an ideal, or an element of $F_\R$, etc.) 
by $\N_{F/\Q}$, and often abbreviate it to $\N$. 
The same rule will be applied to $\Tr=\Tr_{F/\Q}$.

\section{Lattices and the Eisenstein series}
In this section, we develop a theory of Eisenstein series 
over a fixed number field $F$. 
We denote the number of real (resp.\ complex) places of $F$ 
by $r_1$ (resp.\ $r_2$). 

\subsection{Lattices in $\D_F$}\label{subsec:Lattices}
For each infinite place $v\in S_F$, let $\D_v$ be 
the quadratic division algebra over $F_v$. Thus $\D_v$ is isomorphic to 
$\C$ or $\H$ according to whether $v$ is real or complex. 
Choosing such an isomorphism, we define $j_v\in\D_v$ to be the element 
corresponding to $i\in\C$ or $j\in\H$. 
Then we have that $\D_v=F_v\oplus F_vj_v$ for any $v$. 
We also define the Haar measure on $\D_v$ to be the Lebesgue measure 
(resp.\ $4$ times the Lebesgue measure) for $v\in S_F^1$ 
(resp.\ $v\in S_F^2$), so that the above direct sum decomposition 
preserves the measure. 

Next, let us put $\D_F:=\prod_{v\in S_F}\D_v$ and $j_F:=(j_v)_v\in\D_F$. 
Then $\D_F=F_\R\oplus F_\R j_F$ becomes naturally 
a quadratic algebra over $F_\R$. If we write an element $z$ of $\D_F$ 
as $z=x+yj_F$, we regard $x$ and $y$ as elements of $F_\R$. 

Moreover, for $z=(z_v)_v\in\D_F$, we set 
\[\lVert z\rVert_{\D_F}=\prod_{v\in S_F^1}\lvert z_v\rvert\cdot
\prod_{v\in S_F^2}\lvert z_v\rvert^2, \]
where $\lvert\,\cdot\,\rvert$ is the usual absolute value in $\C$ or $\H$. 
Note that $\lVert x\rVert=\absN{x}$ for $x\in F_\R$. 

\begin{defn}\upshape
We call a discrete and cocompact $O_F$-submodule of $\D_F$ 
an \textit{$O_F$-lattice} in $\D_F$. 
For such a lattice $\Lambda$, we denote by $V(\Lambda)$ 
the volume of the quotient $\D_F/\Lambda$ 
(with respect to the product Haar measure on $\D_F$). 
\end{defn}

\begin{lem}\label{lem:V(Lambda)}
Let $\Lambda\subset\D_F$ be an $O_F$-lattice. 
\begin{enumerate}
\item For any $z\in\D_F^\times$, 
we have $V(z\Lambda)=V(\Lambda z)=\lVert z\rVert^2V(\Lambda)$. 
\item There exist elements $\omega_1,\,\omega_2\in\D_F$ and 
a fractional ideal $\fra$ of $F$ such that 
\[\Lambda=\fra\omega_1+O_F\omega_2. \]
\item Let $\fra$ and $\frb$ be fractional ideals of $F$, 
and $z$ an element of $\D_F$ of the form $z=x+yj_F$ where $x\in F_\R$ and 
$y\in F_\R^\times$. Then $\Lambda=\fra z+\frb$ is an $O_F$-lattice, and 
\[V(\Lambda)=d_F\N(\fra)\N(\frb)\absN{y}.\]
\end{enumerate}
\end{lem}
\begin{proof}
The first assertion is clear from the definition. 
(2) is a special case of the structure theorem 
for finitely generated torsion-free modules over Dedekind domains 
(Bourbaki \cite[Chap.~7, \S4, Proposition 24]{Bourbaki85}). 

For (3), we consider the map 
\[F_\R\oplus F_\R\ni(\alpha,\beta)\longmapsto\alpha z+\beta\in\D_F. \]
This is an isomorphism of $\R$-vector spaces, and multiplies the volume 
by $\absN{y}$. On the other hand, we see that 
the volume of $F_\R/\fra$ is $\N(\fra)\sqrt{d_F}$ 
and similar for $\frb$. Then the claim follows. 
\end{proof}

\subsection{The Eisenstein series}\label{subsec:Eisenstein}
\begin{defn}\upshape
An $O_F$-lattice $\Lambda\subset\D_F$ is said to be 
\textit{non-degenerate} if there is no nonzero element $\lambda\in\Lambda$ 
satisfying $\lVert\lambda\rVert=0$. 
For such $\Lambda$, we define 
the \textit{Eisenstein series} $E(\Lambda,s)$ by 
\[E(\Lambda,s):=\sumpr_{\lambda\in\Lambda/U_F}
\frac{V(\Lambda)^s}{\lVert\lambda\rVert^{2s}}\qquad(\re(s)>1),\]
where the prime means that the sum is taken for nonzero $\lambda$. 
\end{defn}

$E(\Lambda,s)$ has the `modularity': 

\begin{lem}\label{lem:Modular}
For $z\in\D_F^\times$, we have $E(z\Lambda,s)=E(\Lambda z,s)=E(\Lambda,s)$. 
\end{lem}
\begin{proof}
This follows from the definition, and Lemma \ref{lem:V(Lambda)} (1). 
\end{proof}

\begin{rem}\upshape\label{rem:Equivalent}
By the above lemma and Lemma \ref{lem:V(Lambda)} (2), (3), 
it is sufficient to consider 
\[E(z,\fra,\frb,s)=\sumpr_{(\mu,\nu)\in(\fra\oplus\frb)/U_F}
\frac{\Bigl(d_F\N(\fra)\N(\frb)\bigl|\N(y)\bigr|\Bigr)^s}
{\lVert\mu z+\nu\rVert^{2s}},\]
for $z=x+yj_F$, $x\in F_\R$ and $y\in F_\R^\times$. 
Moreover, this series can be regarded as a function of 
$\bigl(x_v+\lvert y_v\rvert j_v\bigr)_v$ (and $s$), i.e., 
a function on the product of $r_1$ copies of the upper half planes and
$r_2$ copies of the hyperbolic $3$-space. 
Hence our notion of Eisenstein series is essentially equivalent to 
the more traditional one (see Asai \cite{Asai70}, 
Jorgenson-Lang \cite{Jorgenson-Lang99} or Yoshida \cite{Yoshida03}). 
\end{rem}

To consider the inverse Mellin transform of the Eisenstein series, 
we need some definitions. 

For each infinite place $v\in S_F$, 
put $n_v=1$ or $2$ according to whether $v$ is real or complex. 
We also set $T_F:=F_\R^\times$, and choose the Haar measure 
$\dtimes t$ on it to be $\prod_vdt_v/|t_v|^{n_v}$. 
Then we define 
\begin{gather*}
f(z):=\prod_{v\in S_F}\exp\bigl(-n_v\pi\,\lvert z_v\rvert^2\bigr)
\qquad(z\in\D_F),\\
\Gamma_F(s):=\int_{T_F}f(t)\bigl|\N(t)\bigr|^s\dtimes t,\qquad
\widehat{E}(\Lambda,s):=\Gamma_F(2s)E(\Lambda,s). 
\end{gather*}

\begin{prop}\label{prop:Mellin}
\[\widehat{E}(\Lambda,s)
=V(\Lambda)^s\int_{T_F/U_F}\sumpr_{\lambda\in\Lambda}f(t\lambda)
\bigl|\N(t)\bigr|^{2s}\dtimes t. \]
\end{prop}
\begin{proof}
Putting $\lambda=\mu+\nu j_F$, we have 
\begin{align*}
\Gamma_F(2s)\lVert\lambda\rVert^{-2s}
&=\prod_{v\in S_F}\int_{F_v^\times}
\exp\bigl(-n_v\pi\,\lvert t_v\rvert^2\bigr)\lvert t_v\rvert^{2n_vs}
\bigl(\lvert\mu_v\rvert^2+\lvert\nu_v\rvert^2\bigr)^{-n_vs}
\frac{dt_v}{\lvert t_v\rvert^{n_v}}\\
&=\prod_{v\in S_F}\int_{F_v^\times}
\exp\Bigl(-n_v\pi\,\lvert t_v\rvert^2
\bigl(\lvert\mu_v\rvert^2+\lvert\nu_v\rvert^2\bigr)\Bigr)
\lvert t_v\rvert^{2n_vs}\frac{dt_v}{\lvert t_v\rvert^{n_v}}\\
&=\int_{T_F}f(t\lambda)\bigl|\N(t)\bigr|^{2s}\dtimes t. 
\end{align*}
Hence, by taking the sum and transforming as 
\[\sumpr_{\lambda\in\Lambda/U_F}\int_{T_F}
=\int_{T_F}\sumpr_{\lambda\in\Lambda/U_F}
=\int_{T_F/U_F}\sumpr_{\lambda\in\Lambda}, \]
we obtain the result. 
\end{proof}

\subsection{The functional equation}
\begin{defn}\upshape
Let $\psi_{\D_F}\colon\D_F\longrightarrow\C^\times_1$ be the character 
defined by 
\[\psi_{\D_F}(x+yj_F):=\exp\bigl(2\pi i\Tr_{F/\Q}(x)\bigr).\]
Then, for an $O_F$-lattice $\Lambda\subset\D_F$, 
we define the \textit{dual lattice} $\Lambda^*$ by 
\[\Lambda^*:=\bigl\{\lambda^*\in\D_F\bigm|\psi_{\D_F}(\lambda\lambda^*)=1
\ (\forall \lambda\in\Lambda)\bigr\}.\]
\end{defn}

\begin{prop}\label{prop:Poisson1}
For a non-degenerate $O_F$-lattice $\Lambda\subset\D_F$, set 
\[\Theta(t,\Lambda):=\sum_{\lambda\in\Lambda}f(t\lambda)\qquad(t\in T_F). \]
Then we have 
\[\Theta(t,\Lambda)=V(\Lambda)^{-1}\bigl|\N(t)\bigr|^{-2}
\Theta(t^{-1},\Lambda^*). \]
In particular, $V(\Lambda^*)=V(\Lambda)^{-1}$ holds. 
\end{prop}
\begin{proof}
If we put $f_t(z)=f(tz)$, its Fourier transform is given by 
\[\hat{f}_t(w):=\int_{\D_F}f_t(z)\overline{\psi_{\D_F}(zw)}dz
=\bigl|\N(t)\bigr|^{-2}f_{t^{-1}}(w). \]
Therefore we can use the Poisson summation formula to get 
\[\sum_{\lambda\in\Lambda}f_t(\lambda)=V(\Lambda)^{-1}\bigl|\N(t)\bigr|^{-2}
\sum_{\lambda^*\in\Lambda^*}f_{t^{-1}}(\lambda^*)\]
as desired. 
The last assertion is shown by applying this formula to $\Lambda^*$. 
\end{proof}

\begin{thm}\label{thm:FE}
$\widehat{E}(\Lambda,s)$ can be continued meromorphically to 
the whole $s$-plane, and satisfies the functional equation 
\[\widehat{E}(\Lambda,s)=\widehat{E}(\Lambda^*,1-s).\]
\end{thm}
\begin{proof}
Proposition \ref{prop:Mellin} says that 
\[\widehat{E}(\Lambda,s)=V(\Lambda)^s\int_{T_F/U_F}
\bigl(\Theta(t,\Lambda)-1\bigr)\bigl|\N(t)\bigr|^{2s}\dtimes t.\]
We decompose the integral as $\int_{T_F/U_F}
=\int_{\lvert\N(t)\rvert\geq 1}+\int_{\lvert\N(t)\rvert\leq 1}$. 
Then the former integral converges for every $s\in\C$. On the other hand, 
after inverting the variable $t$, the latter gives 
\begin{align*}
&V(\Lambda)^s\int_{\lvert\N(t)\rvert\geq 1}
\bigl(\Theta(t^{-1},\Lambda)-1\bigr)\bigl|\N(t)\bigr|^{-2s}\dtimes t\\
&=V(\Lambda^*)^{1-s}\int_{\lvert\N(t)\rvert\geq 1}
\bigl(\Theta(t,\Lambda^*)-1\bigr)\bigl|\N(t)\bigr|^{2-2s}\dtimes t
+R,
\end{align*}
where 
\[R=V(\Lambda^*)^{1-s}\int_{\lvert\N(t)\rvert\geq 1}
\bigl|\N(t)\bigr|^{2-2s}\dtimes t
-V(\Lambda)^s\int_{\lvert\N(t)\rvert\geq 1}\bigl|\N(t)\bigr|^{-2s}\dtimes t.\]
Moreover, we see that 
\[\int_{\lvert\N(t)\rvert\geq 1}\bigl|\N(t)\bigr|^{-2s}\dtimes t
=C\int_1^\infty t^{-2s}\frac{dt}{t}=\frac{C}{2s},\]
denoting by $C$ the volume of 
$\bigl\{t\in T_F/U_F\bigm|\lvert\N(t)\rvert=1\bigr\}$ with respect to a 
suitably normalized measure. 
(In fact, we know the value of $C=C_F$. See \ref{subsec:EKLF} below.) 
Hence we obtain the analytic continuation of $R$, and 
the functional equation follows immediately. 
\end{proof}

\subsection{The Fourier expansion}
Here we give the `Fourier expansion' of the Eisenstein series. 
More precisely, we consider $O_F$-lattices of the form 
$\Lambda=\fra z+\frb$, where $\fra$ and $\frb$ are fractional ideals 
and $z=x+yj_F\in\D_F$. Then $\widehat{E}(\Lambda,s)$ is invariant 
under the transforms $x\longmapsto x+\beta$ for all $\beta\in\frb$, 
and hence has the Fourier expansion. 

Let us begin with some preparations. 

\begin{defn}\upshape
\begin{enumerate}
\item To a fractional ideal $\fra$ of $F$, we attach the 
(\textit{completed}) \textit{zeta function} 
\[\zeta_F(s,\fra):=
\N(\fra)^s\sumpr_{\alpha\in\fra/U_F}\bigl|\N(\alpha)\bigr|^{-s},\qquad
\xi_F(s,\fra):=d_F^{s/2}\Gamma_F(s)\zeta_F(s,\fra). \]
\item For $a,b\in F_\R^\times$, we define 
\begin{align*}
B_F(a,b,s)&:=
\int_{T_F}f_t(a)f_{t^{-1}}(b)\bigl|\N(t)\bigr|^{2s}\dtimes t\\
&=(2\pi)^{r_2}\bigl|\N(b/a)\bigr|^s\prod_{v\in S_F}
K_{n_vs}\bigl(n_v\pi\lvert a_vb_v\rvert\bigr). 
\end{align*}
Here $K_s(x)=\int_0^\infty e^{-x(u+u^{-1})}u^{s-1}du$ is the Bessel function. 
\end{enumerate}
\end{defn}

\begin{rem}\upshape
$\zeta_F(s,\fra)$ is equal to the partial zeta function 
\[\zeta_F(s,A)=\sum_{\frb\in A,\,\frb\subset O_F}\N(\frb)^{-s},\]
where $A$ is the wide ideal class containing $\fra^{-1}$. 
\end{rem}

\begin{prop}\label{prop:ZetaMellin}
We have 
\[\xi_F(s,\fra)=V(\fra)^s
\int_{T_F/U_F}\sumpr_{\alpha\in\fra}f(t\alpha)
\bigl|\N(t)\bigr|^s\dtimes t. \]
Here $V(\fra)=d_F^{1/2}\N(\fra)$ denotes the volume of $F_\R/\fra$. 
\end{prop}
\begin{proof}
This can be shown in the same way as Proposition \ref{prop:Mellin}. 
\end{proof}

\begin{lem}\label{lem:Poisson2}
For $z=x+yj_F\in\D_F$, $t\in T_F$ and a fractional ideal $\frb$ of $F$, 
we have 
\[\sum_{\beta\in\frb}f_t(z+\beta)
=V(\frb)^{-1}\bigl|\N(t)\bigr|^{-1}
\sum_{\beta^*\in\frb^*}e^{2\pi i\Tr(x\beta^*)}f_t(y)f_{t^{-1}}(\beta^*).\]
\end{lem}
\begin{proof}
The Fourier transform of the function 
\[f_{t,z}(u):=f\bigl(t(z+u)\bigr)=f_t(y)f_t(x+u)\]
 on $F_\R$ is given by 
\[\hat{f}_{t,z}(v):=\int_{F_\R}f_{t,z}(u)e^{-2\pi i\Tr(uv)}du
=e^{2\pi i\Tr(xv)}\bigl|\N(t)\bigr|^{-1}f_t(y)f_{t^{-1}}(v).\]
Hence the claim is a consequence of the Poisson summation formula. 
\end{proof}

\begin{thm}\label{thm:FourierExp}
Let $\fra$ and $\frb$ be fractional ideals of $F$, and 
$z=x+yj_F$ an element of $\D_F$ with $y\in F_\R^\times$. 
Then, for $\Lambda=\fra z+\frb$, we have 
\begin{align*}
\widehat{E}&(\Lambda,s)\\
=&\biggl(\frac{\N(\fra)}{\N(\frb)}\absN{y}\biggr)^s\xi_F(2s,\frb)
+\biggl(\frac{\N(\fra)}{\N(\frb)}\absN{y}\biggr)^{1-s}
\xi_F(2s-1,\fra)\\
&+V(\fra)^sV(\frb)^{s-1}\absN{y}^s
\sum_{(\alpha,\beta^*)}
e^{2\pi i\Tr(x\alpha\beta^*)}B_F\biggl(\alpha y,\beta^*,s-\frac{1}{2}\biggr). 
\end{align*}
Here $(\alpha,\beta^*)\in
\bigl(\fra\setminus\{0\}\bigr)\times\bigl(\frb^*\setminus\{0\}\bigr)$ 
runs through a system of representatives with respect to 
the equivalence relation defined by 
\[(\alpha,\beta^*)\sim (\alpha\eps,\beta^*\eps^{-1})\qquad
(\forall\eps\in U_F). \]
\end{thm}
\begin{proof}
First, note that 
\[\sumpr_{\lambda\in\Lambda}f_t(\lambda)
=\sumpr_{\beta\in\frb}f_t(\beta)
+\sumpr_{\alpha\in\fra}\sum_{\beta\in\frb}f_t(\alpha z+\beta). \]
Furthermore, by Lemma \ref{lem:Poisson2}, we have 
\begin{align*}
\sumpr_{\alpha\in\fra}\sum_{\beta\in\frb}f_t(\alpha z+\beta)
=&V(\frb)^{-1}\bigl|\N(t)\bigr|^{-1}\sumpr_{\alpha\in\fra}
\sum_{\beta^*\in\frb^*}e^{2\pi i\Tr(x\alpha\beta^*)}
f_t(y\alpha)f_{t^{-1}}(\beta^*)\\
=&V(\frb)^{-1}\bigl|\N(t)\bigr|^{-1}\sumpr_{\alpha\in\fra}
\sumpr_{\beta^*\in\frb^*}e^{2\pi i\Tr(x\alpha\beta^*)}
f_t(y\alpha)f_{t^{-1}}(\beta^*)\\
&+V(\frb)^{-1}\bigl|\N(t)\bigr|^{-1}\sumpr_{\alpha\in\fra}f_t(y\alpha). 
\end{align*}
Therefore, the theorem is deduced from Proposition \ref{prop:Mellin}, 
Lemma \ref{lem:V(Lambda)} (3), Proposition \ref{prop:ZetaMellin}, 
and 
\begin{align*}
&\int_{T_F/U_F}\sumpr_{\alpha\in\fra}\sumpr_{\beta^*\in\frb^*}
e^{2\pi i\Tr(x\alpha\beta^*)}f_t(y\alpha)f_{t^{-1}}(\beta^*)\absN{t}^{2s-1}
\dtimes t\\
&=\sum_{(\alpha,\beta^*)}e^{2\pi i\Tr(x\alpha\beta^*)}
\int_{T_F}f_t(y\alpha)f_{t^{-1}}(\beta^*)\absN{t}^{2s-1}
\dtimes t. 
\end{align*}
The last follows from the identity 
\[f_t(y\alpha\eps)f_{t^{-1}}(\beta^*\eps^{-1})
=f_{t\eps}(y\alpha)f_{(t\eps)^{-1}}(\beta^*)\]
for each $\eps\in U_F$. 
\end{proof}

\subsection{The Kronecker limit formula}\label{subsec:EKLF}
For a meromorphic function $\varphi(s)$ around $s=\alpha$, 
we denote by $\CT_{s=\alpha}\varphi(s)$ the constant term 
in the Laurent expansion at $s=\alpha$, 
while $\Res_{s=\alpha}\varphi(s)$ means the residue. 
As an application of Theorem \ref{thm:FourierExp}, 
we give a limit formula of Kronecker's type 
which expresses $\CT_{s=1}\widehat{E}(\Lambda,s)$. 

Let us denote the regulator of $F$ by $R_F$, and 
the number of roots of unity in $F$ by $w_F$. Then we put 
\[C_F:=\frac{2^{r_1}(2\pi)^{r_2}R_F}{w_F}. \]

\begin{thm}\label{thm:EKLF}
Let $\Lambda=\fra z+\frb$ be as in Theorem \ref{thm:FourierExp}. 
Then 
\begin{gather*}
\Res_{s=1}\widehat{E}(\Lambda,s)=\frac{C_F}{2},\\
\CT_{s=1}\widehat{E}(\Lambda,s)=\CT_{s=1}\xi_F(s,\fra)
+\frac{C_F}{2}\Biggl(h_F(z,\fra,\frb)-
\log\biggl(\frac{\N(\fra)}{\N(\frb)}\absN{y}\biggr)\Biggr),
\end{gather*}
where the function $h_F$ is defined by 
\begin{align*}
\frac{C_F}{2}\,h_F(z,\fra,\frb)=&\ 
\frac{\N(\fra)}{\N(\frb)}\absN{y}\xi_F(2,\frb)\\
&+V(\fra)\absN{y}\sum_{(\alpha,\beta^*)}
e^{2\pi i\Tr(x\alpha\beta^*)}B_F\biggl(\alpha y,\beta^*,\frac{1}{2}\biggr).
\end{align*}
\end{thm}
\begin{proof}
It is well-known that 
\[\Res_{s=1}\zeta_F(s,\fra)=\frac{C_F}{d_F^{1/2}},\]
or equivalently 
\[\Res_{s=1}\xi_F(s,\fra)=C_F. \]
Thus the claim follows from Theorem \ref{thm:FourierExp}. 
\end{proof}

The function $h_F(z,\fra,\frb)$ has a modular property. 

\begin{cor}\label{cor:hModular}
For $\begin{pmatrix}a&b\\c&d\end{pmatrix}\in\mathrm{GL}_2(O_F)$ 
such that $b\in\fra\frb^{-1}$ and $c\in\fra^{-1}\frb$, we have 
\[h_F\bigl((az+b)(cz+d)^{-1},\fra,\frb\bigr)
=h_F(z,\fra,\frb)-2\log\lVert cz+d\rVert. \]
\end{cor}
\begin{proof}
Put $z'=x'+y'j_F=(az+b)(cz+d)^{-1}$. 
The conditions on $b$ and $c$ ensures that 
\[\fra z+\frb=\fra\,(az+b)+\frb\,(cz+d)=(\fra z'+\frb)\,(cz+d), \]
and hence 
\[\widehat{E}(\fra z+\frb,s)=\widehat{E}(\fra z'+\frb,s). \]
Then we see from Theorem \ref{thm:EKLF} that 
\[h_F(z,\fra,\frb)-\log\absN{y}=h_F(z',\fra,\frb)-\log\absN{y'}. \]
Hence the claim is reduced to the identity $\N(y')=\N(y)/\lVert cz+d\rVert^2$, 
which is easily shown. 
\end{proof}

\begin{rem}\upshape
The function $h_F(z,\fra,\frb)$ is a generalization of 
that studied by Asai \cite{Asai70} and Jorgenson-Lang \cite{Jorgenson-Lang99} 
(up to constant multiplication). 
\end{rem}

\section{Hecke's integral formula for quadratic extensions}
Let $K$ be a quadratic extension of a number field $F$. 
The goal in this section is Theorem \ref{thm:HIF}, which represents 
a zeta function of $K$ as an integral of the Eisenstein series for $F$. 

In the following, we denote the non-trivial $F$-automorphism of $K$ 
by $x\longmapsto x'$. We also use the same notation for the induced maps 
on $K_\R=K\otimes_\Q\R$ or $S_K$. 

\subsection{Hecke's integral formula}\label{subsec:HIF}
First, we show an integral formula in a somewhat abstract setting. 
Let us denote by $\dtimes u$ the quotient Haar measure on $T_K/T_F$, 
i.e.\ the measure which satisfies 
\[\int_{T_K}\phi(t)\dtimes t
=\int_{T_K/T_F}\int_{T_F}\phi(t\tilde{u})\dtimes t\,\dtimes u,\]
where $\phi$ is any integrable function on $T_K$ and 
$\tilde{u}\in T_K$ denotes a lift of $u\in T_K/T_F$. 

We define 
\[g(z):=\prod_{w\in S_K}\exp\bigl(-n_w\pi\lvert z_w\rvert^2\bigr)
\qquad \bigl(n_w:=[K_w:\R]\bigr)\]
for $z=(z_w)_w\in K_\R$. 

\begin{prop}\label{prop:AbstractHIF}
Let $\rho\colon K_\R\longrightarrow\D_F$ be an isomorphism of $F_\R$-modules 
which preserves the Haar measure and satisfies $g(z)=f\bigl(\rho(z)\bigr)$ 
for any $z\in K_\R$. Then, for a fractional ideal $\frA$ of $K$, we have 
\[\xi_K(s,\frA)=\int_{T_K/T_FU_K}
\widehat{E}(\rho(\tilde{u}\frA),s)\dtimes u. \]
Here the lift $\tilde{u}\in T_K$ of $u\in T_K/T_FU_K$ 
is chosen to satisfy $\bigl|\N_{K/\Q}(\tilde{u})\bigr|=1$. 
\end{prop}
\begin{proof}
We apply Proposition \ref{prop:ZetaMellin} and compute as 
\begin{align*}
\xi_K(s,\frA)&=V(\frA)^s
\int_{T_K/U_K}\sumpr_{\alpha\in\frA}g(t\alpha)\bigl|\N_{K/\Q}(t)\bigr|^s
\dtimes t\\
&=V(\frA)^s
\int_{T_K/T_FU_K}\int_{T_F/T_F\cap U_K}\sumpr_{\alpha\in\frA}
f\bigl(\rho(t\tilde{u}\alpha)\bigr)\bigl|\N_{K/\Q}(t\tilde{u})\bigr|^s
\dtimes t\,\dtimes u\\
&=\int_{T_K/T_FU_K}V\bigl(\rho(\tilde{u}\frA)\bigr)^s
\int_{T_F/U_F}\sumpr_{\alpha\in\frA}f\bigl(t\rho(\tilde{u}\alpha)\bigr)
\bigl|\N_{F/\Q}(t)\bigr|^{2s}\dtimes t\,\dtimes u. 
\end{align*}
In view of Proposition \ref{prop:Mellin}, this is the desired equation. 
\end{proof}

To obtain a more concrete formula, we need to construct $\rho$, 
describe the group $T_K/T_FU_K$ with its measure, 
and choose a lift $\tilde{u}$ for each $u$. 

For each $v\in S_F$, choose and fix a place $w=w_v\in S_K$ above $v$. 
Then we define a map $\rho\colon K_\R\longrightarrow\D_F$ by 
\[\bigl(\rho(z)\bigr)_v:=z_w+z_{w'}j_v
=\begin{cases}
z_w+z_{w'}i&(v\in S_F^1,\,w\in S_K^1),\\
(1+i)z_w&(v\in S_F^1,\,w\in S_K^2),\\
z_w+z_{w'}j&(v\in S_F^2,\,w\in S_K^2).\end{cases}\]
It is easy to check that it satisfies the conditions 
in Proposition \ref{prop:AbstractHIF}. 

Next, we set 
\[T_{K/F}:=\bigl\{u\in T_K\bigm|\N_{K/F}(u)=1\bigr\},\quad 
U_{K/F}:=T_{K/F}\cap U_K. \]
Then the map $x\longmapsto x/x'$ induces the isomorphisms 
\[T_K/T_F\overset{\cong}{\longrightarrow}T_{K/F},\qquad 
T_K/T_FU_{K/F}\overset{\cong}{\longrightarrow}T_{K/F}/U_{K/F}^2,\]
where $U_{K/F}^2$ means $\{u^2\mid u\in U_{K/F}\}$. 

Let us study the structure of $T_{K/F}$ by looking at its $v$-component 
\[(T_{K/F})_v=\begin{cases}
\bigl\{\bigl(z,z^{-1}\bigr)\in K_w^\times\times K_{w'}^\times\bigm| 
z\in F_v^\times\bigr\}&(F_v=K_w=K_{w'}),\\
\bigl\{z\in K_w^\times\bigm| \lvert z\rvert=1\bigr\}
&(F_v\subsetneq K_w\cong\C).
\end{cases}\]
for each $v\in S_F$. 
When $F_v=K_w$, we choose the Haar measure on $(T_{K/F})_v$ so that 
the isomorphism 
\[F_v^\times\ni z\longmapsto\bigl(z,z^{-1}\bigr)\in(T_{K/F})_v\]
preserves the measure. 
Moreover, for $u=\bigl(z,z^{-1}\bigr)\in(T_{K/F})_v$, we define 
\[\tilde{u}:=\bigl(z\,\lvert z\rvert^{-1/2},\lvert z\rvert^{-1/2}\bigr)
\in K_w^\times\times K_{w'}^\times. \]
On the other hand, if $F_v\subsetneq K_w$, 
we equip $(T_{K/F})_v\cong\C^\times_1$ the measure of total mass $2\pi$, 
and put $\tilde{u}:=\sqrt{u}\in K_w^\times$, any one of the square roots 
of $u\in(T_{K/F})_v$. 

Now, by taking the product, we obtain the Haar measure $\dtimes u$ 
on $T_{K/F}$ and the map 
\[T_{K/F}\ni u\longmapsto \tilde{u}\in T_K\]
satisfying $\tilde{u}/\tilde{u}'=u$ and $\absN{\tilde{u}}=1$. 

\begin{thm}\label{thm:HIF}
For a fractional ideal $\frA$ of $K$, we have 
\[\xi_K(s,\frA)=\frac{1}{w_{K/F}}
\int_{T_{K/F}/U_{K/F}^2}\widehat{E}\bigl(\rho(\tilde{u}\frA),s\bigr)
\dtimes u, \]
where $w_{K/F}$ denotes the index $[U_K:U_FU_{K/F}]$. 
\end{thm}
\begin{proof}
We have only to check the compatibility of Haar measures 
in the isomorphism $T_K/T_F\cong T_{K/F}$, i.e.\ 
\[\int_{T_K}\phi(t)\dtimes t
=\int_{T_{K/F}}\int_{T_F}\phi(t\tilde{u})\dtimes t\,\dtimes u\]
for integrable functions $\phi$. 
Let us consider componentwise. 

If both $v$ and $w$ are real, we have 
\[F_v^\times=K_w^\times=K_{w'}^\times\cong \{\pm 1\}\times\R^\times_+. \]
Here $\R^\times_+$, the multiplicative group of positive real numbers, 
is equiped with the Haar measure $dt/t$, 
while $\{\pm 1\}$ has the total mass $2$. 
Therefore we may consider the compatibilities of measures in the isomorphisms 
\[\bigl(\{\pm 1\}\times\{\pm 1\}\bigr)/\{\pm 1\}\cong\{\pm 1\},\qquad
\bigl(\R^\times_+\times\R^\times_+\bigr)/\R^\times_+\cong\R^\times_+\]
separately. Then the problem in the former is trivial, 
while that in the latter is reduced to the elementary formula 
\[\int_0^\infty\int_0^\infty \phi(x,y)dx\,dy=
\int_0^\infty\int_0^\infty \phi\bigl(tu^{1/2},tu^{-1/2}\bigr)dt\,du. \]
The case in which $v$ and $w$ are complex can be 
treated in the same way, using 
$F_v^\times\cong \C^\times_1\times\R^\times_+$, 
where $\C^\times_1$ has the total mass $4\pi$ (recall that the Haar measure 
on $F_v$ and $K_w$ are twice the Lebesgue measure). 
Finally, if $v$ is real and $w$ is complex, we have 
\[(T_{K/F})_v=K_w^\times/F_v^\times
\cong\bigl(\C^\times_1\times\R^\times_+\bigr)\bigm/
\bigl(\{\pm 1\}\times\R^\times_+\bigr)
\cong\C^\times_1/\{\pm 1\},\]
and both sides have the total mass $2\pi$. This completes the proof. 
\end{proof}

\begin{rem}\upshape
Instead of $\rho$, we can use the map $\rho^*$ defined by 
\[\bigl(\rho^*(z)\bigr)_v:=z_w-j_vz_{w'}
=\begin{cases}
z_w-z_{w'}i&(v\in S_F^1,\,w\in S_K^1),\\
(1-i)z_w&(v\in S_F^1,\,w\in S_K^2),\\
z_w-\overline{z_{w'}}j&(v\in S_F^2,\,w\in S_K^2).\end{cases}\]
Then the dual lattice of $\rho(\tilde{u}\frA)$ is 
$\rho^*(\tilde{u}^{-1}\frA^*)$. 
Hence the functional equation for $\widehat{E}$ (Theorem \ref{thm:FE}) 
leads to the functional equation $\xi_K(s,\frA)=\xi_K(1-s,\frA^*)$. 
\end{rem}

\subsection{Application to the Kronecker limit formula}\label{subsec:zetaKLF}
Combining Theorem \ref{thm:HIF} and Theorem \ref{thm:EKLF}, 
we obtain a `relative' Kronecker limit formula, 
which represents a relation between $\CT_{s=1}\xi_K(s,\frA)$ and 
$\CT_{s=1}\xi_F(s,\fra)$. 

\begin{thm}\label{thm:zetaKLF}
Let $\frA\subset K$ be a fractional ideal of the form 
$\frA=\fra z+\frb$ where $\fra$ and $\frb$ are fractional ideals of $F$. 
Then we have 
\begin{align*}
\frac{\CT_{s=1}\xi_K(s,\frA)}{C_K}=&\ 
2\frac{\CT_{s=1}\xi_F(s,\fra)}{C_F}-\log\frac{\N(\fra)}{\N(\frb)}\\
&+\frac{C_F}{2w_{K/F}C_K}\int_{T_{K/F}/U_{K/F}^2}\Bigl(
h_F(z_{\tilde{u}},\fra,\frb)-\log\absN{y_{\tilde{u}}}
\Bigr)\dtimes u. 
\end{align*}
Here we put 
\[z_{\tilde{u}}=x_{\tilde{u}}+y_{\tilde{u}}j_F
:=\rho(\tilde{u}z)\rho(\tilde{u})^{-1}.\]
\end{thm}
\begin{proof}
Since 
\[\rho(\tilde{u}\frA)=\fra\,\rho(\tilde{u}z)+\frb\,\rho(\tilde{u})
=(\fra z_{\tilde{u}}+\frb)\,\rho(\tilde{u}), \]
Theorem \ref{thm:HIF} tells us that 
\[\xi_K(s,\frA)=\frac{1}{w_{K/F}}
\int_{T_{K/F}/U_{K/F}^2}\widehat{E}(\fra z_{\tilde{u}}+\frb,s)\dtimes u. \]
By comparing the residues at $s=1$, we obtain 
\[C_K=\frac{C_F}{2w_{K/F}}\int_{T_{K/F}/U_{K/F}^2}\dtimes u. \]
Thus the claimed formula follows from Theorem \ref{thm:EKLF}. 
\end{proof}



\begin{thebibliography}{99}
\bibitem{Asai70} T.~Asai, 
	On a certain function analogous to $\log\lvert\eta(z)\rvert$, 
	Nagoya Math.\ J., 40 (1970), 193--211. 
\bibitem{Bourbaki85} N.~Bourbaki, 
	Alg\`ebre Commutative, 
	Masson, 1985. 
\bibitem{Hecke17} E.~Hecke, 
	\"Uber die Kroneckersche Grenzformel f\"ur reelle quadratische K\"orper 
	und die Klassenzahl relative-Abelscher K\"orper, 
	Mathematische Werke, 198--207. 
\bibitem{Jorgenson-Lang99} J.~Jorgenson, S.~Lang, 
	Hilbert-Asai Eisenstein series, regularized products, 
	and heat kernels, 
	Nagoya Math.\ J., 153 (1999), 155--188. 
\bibitem{Konno89} S.~Konno, 
	Eisenstein series in hyperbolic $3$-space and Kronecker limit formula 
	for biquadratic field, 
	Nagoya Math.\ J., 113 (1989), 129--146. 
\bibitem{Manin04} Y.~I.~Manin, 
	Real multiplication and noncommutative geometry (ein Alterstraum), 
	\textit{The Legacy of Niels Henrik Abel}, 685--727, Springer, Berlin, 2004.
\bibitem{Meyer57} C.~Meyer, 
	Die Berechnung der Klassenzahl Abelscher K\"orper 
	\"uber quadratischen Zahlk\"orpern, 
	Akademie-Verlag, Berlin, 1957. 
\bibitem{Siegel80} C.~L.~Siegel, 
	Advanced Analytic Number Theory, 
	Tata Institute of Fundamental Research, Bombay, 1980. 
\bibitem{Yoshida03} H.~Yoshida, 
	Absolute CM-Periods, 
	Mathematical Surveys and Monographs, Vol.~106, 
	American Mathematical Society, 2003. 
\end{thebibliography}
\end{document}